\documentclass[conference]{IEEEtran}

\ifCLASSINFOpdf

\else

\fi
\usepackage{graphicx,algorithm2e}
\DeclareGraphicsExtensions{.pdf,.png,.jpg}
\begin{document}

\title{Simulated Tornado Optimization}
% when the number of particles are such low then the turbulance would not be better than normal (Gaussian)
\author{\IEEEauthorblockN{S. Hossein Hosseini$^\ast$, Tohid Nouri$^\dag$, Afshin Ebrahimi$^\ast$, and S. Ali Hosseini$^\ddag$}
\IEEEauthorblockA{$^\ast$ICT Research Center, Faculty of Electrical Engineering, Sahand University of Technology, Tabriz, Iran\\
$^\dag$Faculty of Civil Engineering, University of Mohaghegh Ardabili, Ardabil, Iran\\
$^\ddag$Faculty of Industrial Engineering, Urmia University of Technology, Urmia, Iran\\
Emails: \{h\_hosseini, aebrahimi\}@sut.ac.ir, t.nouri@uma.ac.ir, syd.ali.hos@gmail.com}}

\IEEEoverridecommandlockouts
%\IEEEpubid{\makebox[\columnwidth]{978-1-5090-5820-4/16/\$31.00~\copyright2016~IEEE\hfill}
%\hspace{\columnsep}\makebox[\columnwidth]{}}
%\IEEEpubid{978-1-5090-5820-4/16/\$31.00~\copyright2016~IEEE\hfill}
\IEEEpubid{}

\maketitle

\begin{abstract}
We propose a swarm-based optimization algorithm inspired by air currents of a tornado. Two main air currents - spiral and updraft - are mimicked. Spiral motion is designed for exploration of new search areas and updraft movements is deployed for exploitation of a promising candidate solution. Assignment of just one search direction to each particle at each iteration, leads to low computational complexity of the proposed algorithm respect to the conventional algorithms. Regardless of the step size parameters, the only parameter of the proposed algorithm, called tornado diameter, can be efficiently adjusted by randomization. Numerical results over six different benchmark cost functions indicate comparable and, in some cases, better  performance of the proposed algorithm respect to some other metaheuristics.
\end{abstract}


\IEEEpeerreviewmaketitle

\section{Introduction}

Metaheuristic or evolutionary algorithms are a family of optimization algorithms inspired by nature, art, or social developments. Despite of their lack in convergence guarantee, these kind of algorithms have been popular in various engineering applications and scientific researches because of their simplicity and efficiency in finding global optimum solution for relatively low-dimension problems. Metaheuristics can be divided to single-solution approaches and multi-solution methods. The well-known examples of the single-solution approach are simulated annealing \cite{SA83} and tabu search \cite{Tabu89} algorithms in which one solution vector is evolved through strategic steps. On the other hand, in multi-solution(agent) methods a number of solutions are appointed and their relations are defined in an evolutionary or optimistic way such that they approach hopefully to a globally optimal solution. Hence, the main difference among the algorithms in the multi-solution approaches is in the kind of relations (motions) defined among the solutions.

Multi-solution approaches are recognized by two prior works on genetic algorithm (GA) \cite{Holland92} and particle swarm optimization (PSO) \cite{JK95}. Other algorithms, for example, differential evolution \cite{RS97}, imperialist competition algorithm (ICA) \cite{EA07}, and teaching-learning-based optimization (TLBO) \cite{TLBO11} are evolved versions of the mentioned two algorithms. They evolve to handle an smart balance between exploration and exploitation. The main differences are in the update rule and selection criteria. Specifically, in the swarm-based approaches the main difference is in the exploration strategy, since a motion toward a globally best solution(s) is a fixed part of these algorithms. The exploration item in PSO is accomplished by assigning a memory to the particles to save the best result of their search. In the ICA, the exploration was attended by appointing several imperialists selected from best solutions. Finally, in the TLBO, the exploration is a part of an educational system which is accomplished by interaction between pairs of the students chosen at random. Our proposed algorithm is another swarm-based optimization method in which the relation and movements among dull air particles of the tornado are mimicked. At the proposed algorithm, the exploration rule is designed inspired by spiral motion of the tornado. Our motivation behind modeling of the tornado motions was its powerful capability in the \emph{reduction} of ambient temperature (up to $15^{\circ}$C) by its mysterious air currents.

In fact, in formulation of a metaheuristic algorithm, there is a wide spectrum between random search and deterministic search where the intensities of exploration and exploitation of an new algorithm is determined. It gives an opportunity to design different algorithms to handle various tradeoffs in facing real-word problems. Two main drawbacks among the metaheuristic algorithms are run-time and parameter adjustability. At this work, both of them, with an emphasis on the former, is addressed. Two search directions is defined to boost the exploration and exploitation capability of the proposed algorithm. Each one is assigned to just one particle. That is despite of conventional methods, e.g. TLBO and PSO algorithms, which apply both main motions to each particle. Consequently, in the proposed algorithm, computational complexity is considerably reduced. Regardless of the step size parameters, the only parameter of the algorithm is the number of particles that should be assigned to each motion. The simulation results indicate that, in most of the problem, this parameter can be floated without significant degradation in the performance. This simple self-adjustment works, since the assignment is done at random.

%In the other point of view, a main difficulty with the mataheuristic algorithms is in their requirement for tuning several parameters in facing with different problems. As another expression, the existence of several parameters in most of the algorithms leads to a diversity in the performance when dealing with different problems. In order to have a general algorithm, hyperheuristic algorithms are developed to manage the parameters or a set of metaheuristic algorithms. That is a higher level program which does the role of a human in the tuning of the parameters or choosing the appropriate algorithm for a specific problem \cite{EB13}. Main drawback of the hyperheuristics is a considerable overhead in evaluation of the low-level metaheuristics or the tuning parameters. That degrades efficiency of the overall algorithm. As a conclusion, development of parameter-free algorithms has special importance because of their generality in solving different problems without need for regular tuning. The proposed algorithm can be regarded as a parameter-free algorithm, in this sense that its motions are defined in a such way that main parameter of the algorithm can be self-adjusted by simply randomization of its value without significant lose of the performance. Simulation results indicate this subject.

\IEEEpubidadjcol

\section{Tornado Morphology}

\begin{figure}[h!]
 \centering
  \includegraphics[scale=0.15]{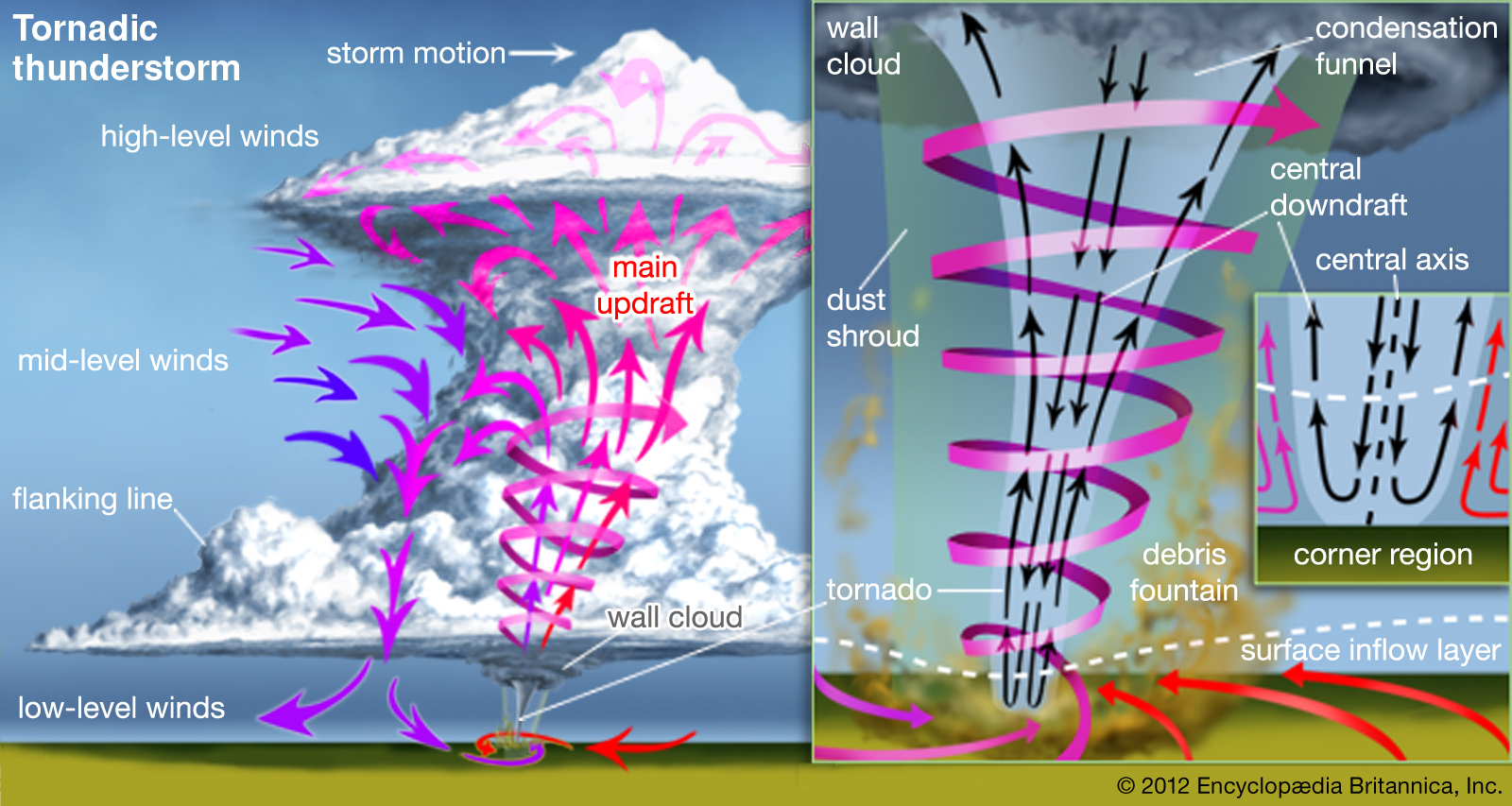}
  \caption{Illustration of air currents in a tornado \cite{thunder16}}
  \label{tormor}
\end{figure}

Tornado is a kind of extreme weather. The strongest wind that could be formed in the earth is related to this phenomenon. When a layer of warm and moist air is located under the cool and dry air, due to the lower density of the warmer one it attempts to climb up to the top of the cold air. Conversely cold air descent to replace the rising air. If this process is occurred quickly, spiral airflow (like a funnel) is formed that is called Tornado \cite{tornado15}. At the beginning of the Tornado the axis of the rotating airflow is horizontal. At the present of wind, the rising warm and moist air tilts the rotating air to vertical. A schematic model of Tornado is shown in Fig. \ref{tormor} \cite{thunder16}. During the Tornado, the warm air blows towards the tornado and rises through the spiral paths. This continues until combination of warm and cold air together and establishment of thermal equilibrium. Two main and specific air currents of the tornado, i.e. updraft and spiral motions, form the basis of our proposed algorithm. Moreover, atmospheric turbulence is an important part of these motions that has significant influence in the proposed simulated tornado optimization (STO) algorithm.

\section{Proposed Algorithm}

\begin{figure}[h!]
 \centering
  \includegraphics[scale=0.24]{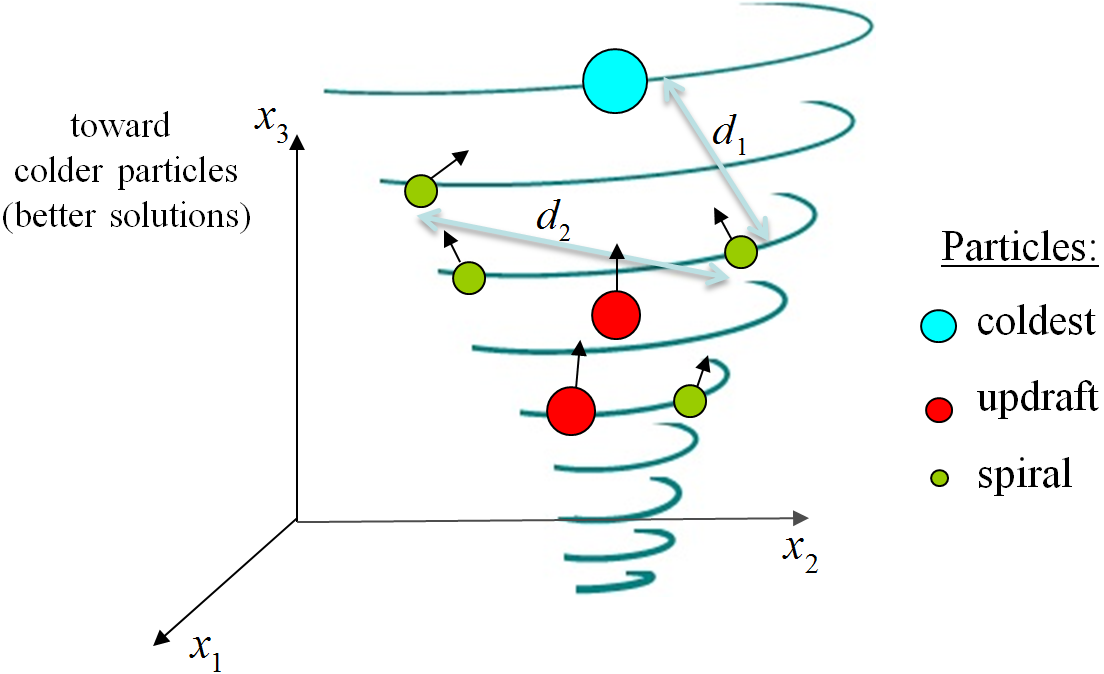}
  \caption{Demonstration of two main movements in the proposed algorithm, regardless of the uncertainty in directions governed by turbulence. ($d_{1}<d_{2}$)}
  \label{STO}
\end{figure}

Let the columns of matrix $X=[x^1~x^2~...~x^k]$ be consist of $k$ solution vectors $x^i=[x_1,~ x_2, ~..., x_N]^T~ (i=1,~...,k$). Each solution includes $N$ decision variables. The solutions are interpreted as the position of $k$ air particles in a tornado. Their corresponding cost function $f(x^i)$ determines the temperature of $i^{th}$ particle at position $x^i$. The particles are divided into three types of coldest ($x^1$), spiral ($i=2,~...,k_1$), and updraft particles ($i=k_1+1,~...,k$), where $1<k_1<k$. Fig. \ref{STO} demonstrates this terminology and illustrates the idea behind the direction of motion for each particle. The coldest particle ($x^1$), in top of all other particles, is best solution at each iteration. For a minimization problem that is a position vector with minimum cost function among the particles. Its replacements, after each iteration, simulates the storm motion (see Fig. \ref{tormor}, top of the left-hand picture). Other positions ($x^i,~i=2,~...,k$) are randomly sorted at the columns of $X$ matrix to participate as a spiral or updraft particle.

Updraft particles (the particles inside the funnel of a real tornado) move directly toward the coldest particle. The following relation formulates this kind of motion:

\begin{equation}\label{asan}
  x^i := x^i+\mu\odot(x^1-x^i)   ~ ~~~\forall~ i=k_1+1,~\ldots,k
\end{equation}

\noindent where $\odot$  indicates an entry-wise multiplication and $\mu$ is the turbulence vector, consist of $N$ step sizes for each decision variable. That is responsible for slight divergences from direct motion toward the coldest particle $x^1$, in order to avoid from

\noindent a pure exploitation. The mentioned formulation is a common rule in the attraction-based algorithms \cite{Yang15}. On the other side, as an special motion - designed to promote the exploration property of the proposed algorithm - a spiral particle at position $x^i~(i=2,~...,k_1)$ moves toward the position of another spiral particle or the coldest one, i.e. $x^j,~j=1,~...,k_1 ~and ~ j\neq i$. This index is determined according to the following criteria:

\begin{equation}\label{asan}
  j=\arg\min ~ \|x^j-x^i\|_2   ~ ~~~s.t. ~~~ f(x^j)<f(x^i)
\end{equation}

\noindent where $\|.\|_2$ denotes Euclidian norm. Hence, each particle moves toward another nearest-better spiral particle or the coldest one (see Fig. \ref{STO}). Corresponding update rule for a spiral motion can be formulated as:

\begin{equation}\label{mes}
  x^i := x^i+\mu\odot(x^j-x^i) ~ ~~~\forall~ i=2,~\ldots,k_1
\end{equation}

\noindent After each iteration, the coldest position is refreshed and spiral-updraft assignments are renewed, randomly.

\subsection{On the Turbulence Model}

The Turbulence vector $\mu$ imposes an uncertainty to both spiral and updraft motion. It is modeled by i.i.d. normal distribution with zero mean and unit variance. This distribution for both kind of updates leads to more generality and acceptable performance respect to the uniform distribution. In addition, that is more general in dealing with different problems respect to the realistic atmospheric turbulence models like log-normal distribution utilized in the communication channel modeling \cite{turbulence10}. Although, that is important to mention, according to our observations, log-normal distribution had shown some promisingly better performance than the normal distribution, at the case of unimodal problems. Nevertheless, for the sake of generality, our simulations are done using normal distribution.

\subsection{Tornado Diameter}

As it is illustrated at the simulations (Fig. \ref{prob}), the number of spiral particles $k_1$ can considerably influence the performance of the proposed algorithm in solving different problems. Its ratio respect to the total number of particles $k$ is called tornado diameter $\alpha=\frac{k_1}{k}$. Each problem has its own optimal value for tornado diameter. As a parameter-free version of the proposed algorithm, we suggest to float this parameter at each iteration to be selected uniformly among the feasible range of $[1,k)$. Simulation results confirm the efficiency of the randomized tornado diameter. The finalized procedure of the proposed algorithm is summarized as follows:

\RestyleAlgo{boxruled}
\LinesNumbered
\begin{algorithm}
\caption{Simulated Tornado Optimization (STO) \label{alg}}
\emph{Generate initial position of $k$ particles by random.}\\
\emph{Evaluate the temperature (cost function) of the particles and appoint the particle with lowest temperature as the coldest particle}\\
  %\item \emph{Assign $k_1$ spiral particles and $k-k_1-1$ updraft particles randomly.
\emph{Select a random integer number for the parameter $k_1$ distributed uniformly between $[1,k)$.}\\
\emph{Randomly assign $k_1$ particles for spiral current and rest for updraft current, then update their positions.}\\
\emph{Repeat 2 to 4 until tornado is vanished (position of the solutions are same).}
\end{algorithm}

\subsection{On the Convergence and Complexity}

The movement toward best solution (coldest particle) in the updraft current is inherently a converging motion \cite{Yang15}. However, that has a risk of falling into a local minimum. The spiral current boosts the exploration capability of the algorithm. Akin to the rings of a chain, the spiral particles are partially connected together and conduce an implicit motion toward the best solution. Hence, the spiral motion itself is a kind of converging attraction. Intuitively, the algorithm still converges after adding the spiral current. However, an analysis of convergence would have insight in choosing the step sizes. It is worth to mention that the randomized assignments of the motions at each iteration is not an essential part of the algorithm. Assignment according to the temperature can also be effective in some problems.

On the other hand, STO has low computational complexity, mainly because the motion of each particle is just in one direction at each iteration. That is despite of the algorithms like TLBO and PSO in which the final direction of movement for each particle at each iteration is determined by the computation of superposition of two or more direction vectors.

\section{Simulation Results}

\begin{figure}[h!]
 \centering
  \includegraphics[scale=0.18]{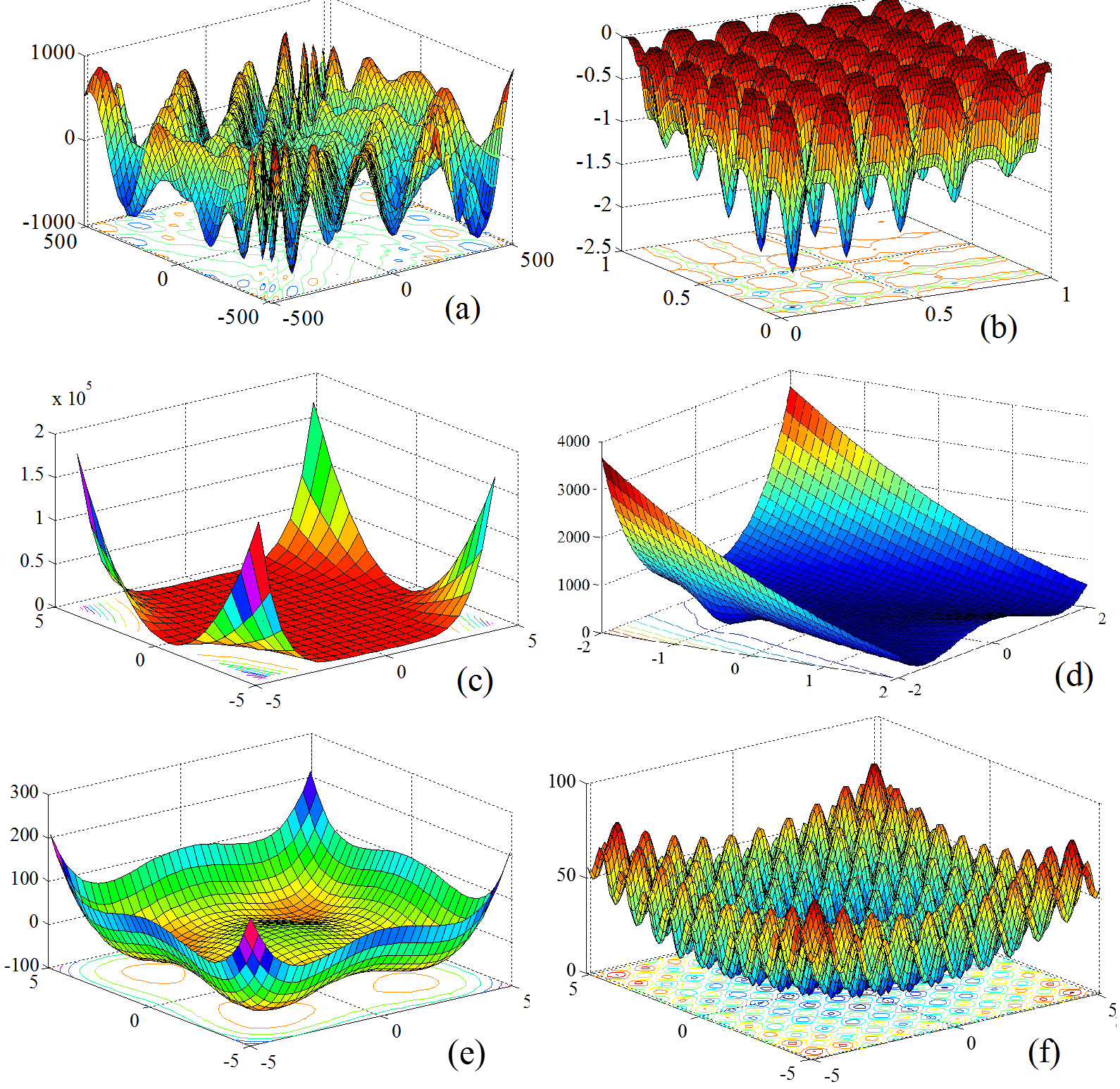}
  \caption{Benchmark cost functions: (a) \emph{Eggholder}, (b) \emph{Ripple}, (c) \emph{Beale}, (d) \emph{Modified Rosenbrock}, (e) \emph{Styblinski-Tang}, (f) \emph{Rastrigin}}
  \label{Rast2}
\end{figure}

At this section, optimization of six different benchmark cost functions are investigated. These functions are shown in Fig. \ref{Rast2} and their formulations, along with domains, optimal solutions, and minimum costs are postponed to appendix. These functions were chosen from \cite{MJ13}. Simulations are conducted through five experiments. Firstly, the trajectory of the particles in the proposed algorithm is demonstrated. Secondly, the effect of tornado diameter on the optimization performance is investigated. Thirdly, influence of the randomization of tornado diameter was considered in a quantitative comparison with a fixed diameter version. Next, fourth experiment includes some qualitative comparisons with with GA, PSO, and TLBO. Finally, the performance of the STO was evaluated on two higher dimension problems.

In all of the experiments the number of particles or population size was fixed on 40 for all of the algorithms. Distortion function is defined as $\frac{\|x^\ast-\hat{x}\|_2}{\|x^\ast\|_2}$, where $x^\ast$ is the real global minimum and $\hat{x}$ is the result of optimization. One run of the optimization algorithms was regarded as a perfect optimization, if the distortion was less than $10^{-4}$. For \emph{Modified Rosenbrock} and \emph{Rastrigin} functions this definition was considered as a cost value less than 36 and $10^{-4}$, respectively. Parameters of the PSO were set according to the constriction coefficients with $\Phi_1=\Phi_2=2.05$ \cite{Kennedy02}. In the GA, 0.8 of the population were selected for crossover and the rest for mutation. Crossover coefficient and mutation probability were tuned on 0.05 and 0.08, respectively. The TLBO is a parameter-free algorithm.

\begin{figure}[h!]
 \centering
  \includegraphics[scale=0.325]{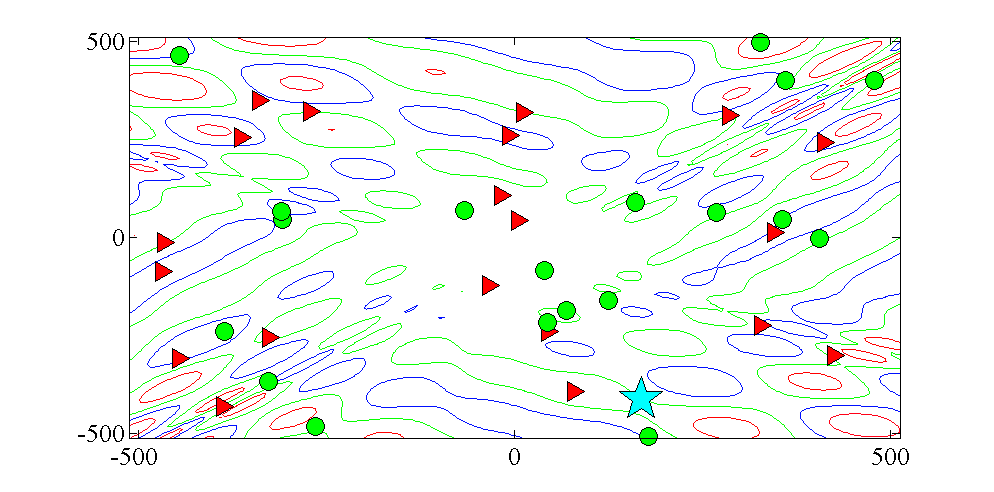}
  \caption{Initial Position of the particles; Star indicator shows the coldest particle and circles and triangles indicate spiral and updraft particles, respectively.}
  \label{traj1}
\end{figure}

\begin{figure}[h!]
 \centering
  \includegraphics[scale=0.325]{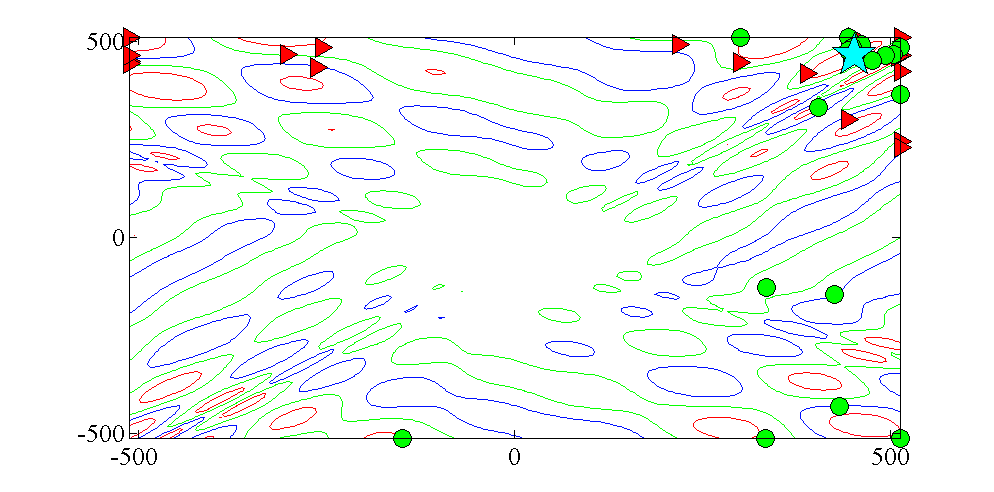}
  \caption{Tornado at iteration 10: Updraft movements toward coldest particle and spiral motion among neighborhoods}
  \label{traj2}
\end{figure}

\begin{figure}[h!]
 \centering
  \includegraphics[scale=0.325]{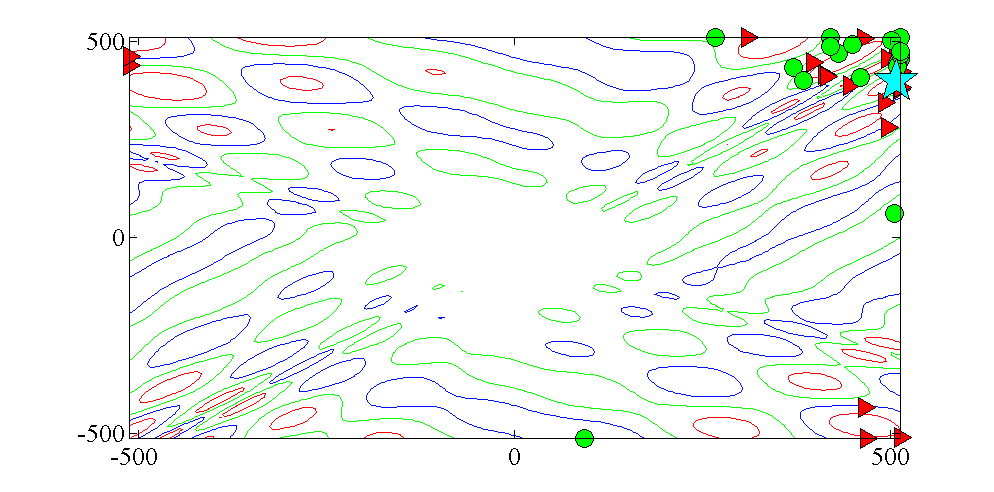}
  \caption{Tornado at iteration 20: coldest particle converges to coldest position (global minimum)}
  \label{traj3}
\end{figure}

\begin{figure}[h!]
 \centering
  \includegraphics[scale=0.325]{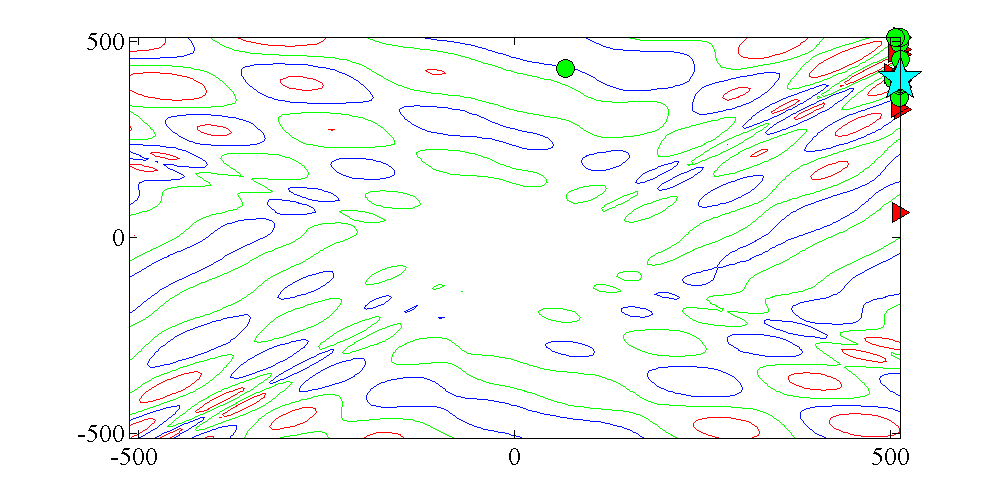}
  \caption{Tornado at iteration 30: All particles approach coldest position and tornado vanishes}
  \label{traj4}
\end{figure}

\begin{figure}[h]
 \centering
  \includegraphics[scale=0.25]{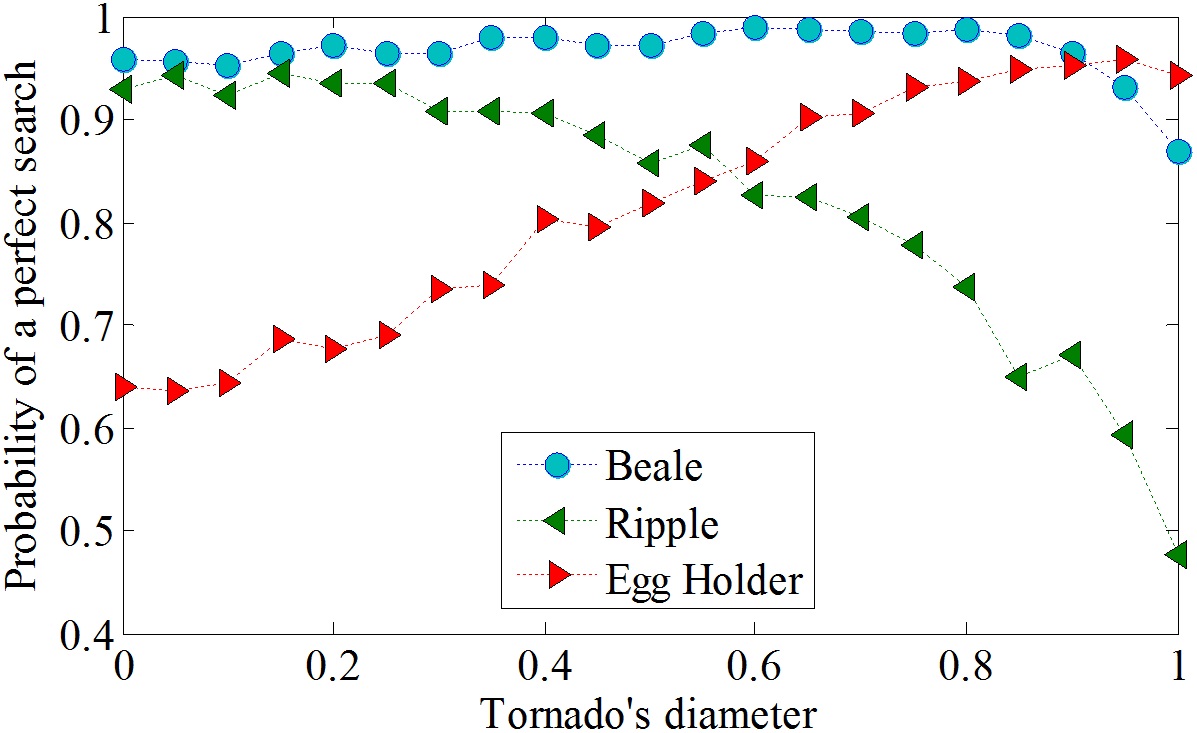}
  \caption{Probability of perfect optimization verse tornado diameter (in each diameter, the number of particles in each air current is fixed)}
  \label{prob}
\end{figure}

At the first experiment, trajectory of the air particles on artificial landscape of \emph{EggHolder} function was demonstrated through Fig. \ref{traj1} to Fig. \ref{traj4}. The number of spiral particles was fixed to be equal with that of updraft particles (including coldest one as a spiral particle). At the second experiment, influence of the parameter of tornado diameter was investigated. For this aim, the algorithm was implemented in whole range of the possible number of the spiral particles $k_1$. Fig. \ref{prob} illustrates the effect of this number in the probability of perfect optimization of three test functions. The success rate was computed over 1000 trials at each values of the $k_1$. As shown, \emph{EggHolder} and \emph{Ripple} functions are two examples of the functions that behave at extreme; the \emph{EggHolder} function was efficiently optimized by high diameter tornados while the optimization of \emph{Ripple} function was more efficient when diameter was small. In between there is \emph{Beale} function which had low sensitivity to the tornado diameter. The proposed algorithm, had best performance in the optimization of \emph{Beale} function when the ratio of the two spiral and updraft particles was in middle points.

\begin{figure}[h!]
 \centering
  \includegraphics[scale=0.29]{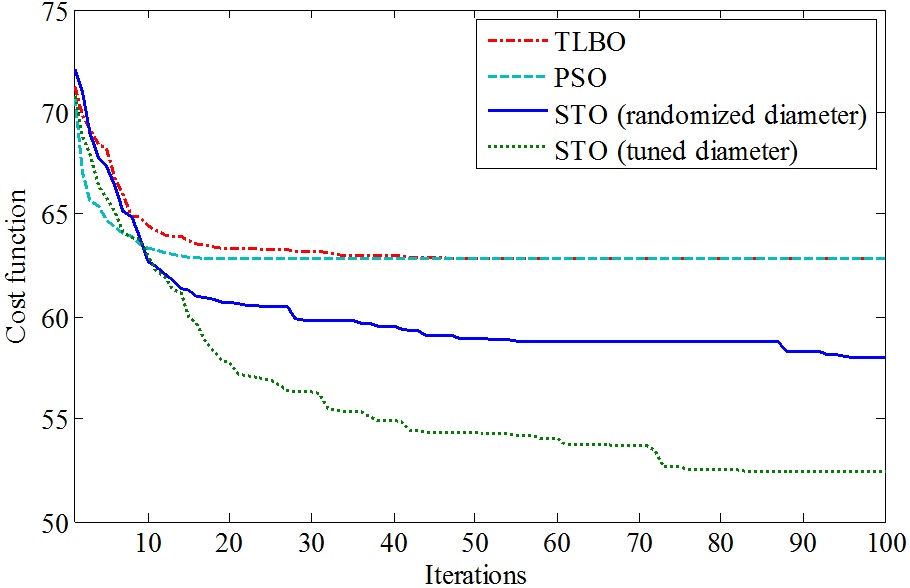}
  \caption{Convergence curves for \emph{Modified Rosenbrock} function}
  \label{sty}
\end{figure}

As mentioned in the previous section, possibility of random assignment of the kind of motion for each particle position, allows to float the number of particles for each air current at each iteration. This randomization is not degrading and in most functions leads to a comparable performance with a tuned version of the STO. Table \ref{tab} compares the probability of perfect optimization at the case of randomized tornado diameter (indicated by STO$^2$) with the best result obtained by tuning of the number of particles assigned to spiral motion (indicated by STO$^1$). As shown, the performances are near (regardless of \emph{Modified Rosenbrock} function). Roughly speaking, possibility of the randomization of parameter $k_1$ at each iteration indicates a promising self-adjustment capability of the proposed algorithm. In addition, the performance of other three algorithms are included at this table. The results are obtained over 1000 trials and all of the algorithms are stopped after 100 iterations. As inferred, the proposed algorithm has comparable quantitative performance in most of the problems and a better one at the case of \emph{Modified Rosenbrock} function. In the rest of experiments, STO in the randomized diameter version was evaluated, except of Fig. \ref{sty}.

\begin{table}[h!]
 \caption{Probability of success in perfect optimization}
\begin{center}
  \begin{tabular}{l c c c c c}
    & \emph{EggHolder} & \emph{Ripple}& \emph{Beale} & \emph{Rosen. Mod.} & \emph{Rast.}(5D) \\
  \hline
  STO\footnote{best performance by tuning}& 0.95 & 0.94 & 0.99& 0.54 ($k_1=35$) & -\\
  STO\footnote{parameter-free version (estimated $k_1$)}& 0.91 & 0.93& 0.98 & 0.40& 0.99 \\
  PSO & 0.24 & 0.26& 0.89 & 0.28& 0.09 \\
  GA & 0.18 & 0.88& 0.47 & 0.18 & 0.46 \\
  TLBO & 0.97 & 0.93& 1 & 0.28& 0.97\\
  \label{tab}

\end{tabular}
\end{center}
\end{table}

\begin{figure}[h!]
 \centering
  \includegraphics[scale=0.30]{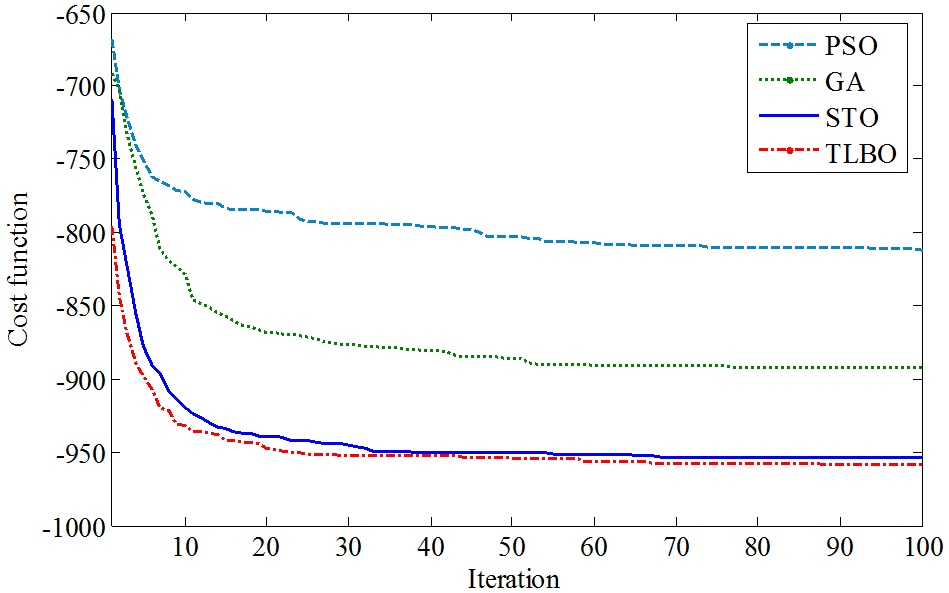}
  \caption{Convergence curves for \emph{EggHolder} function}
  \label{kual1}
\end{figure}

\begin{figure}[h!]
 \centering
  \includegraphics[scale=0.2325]{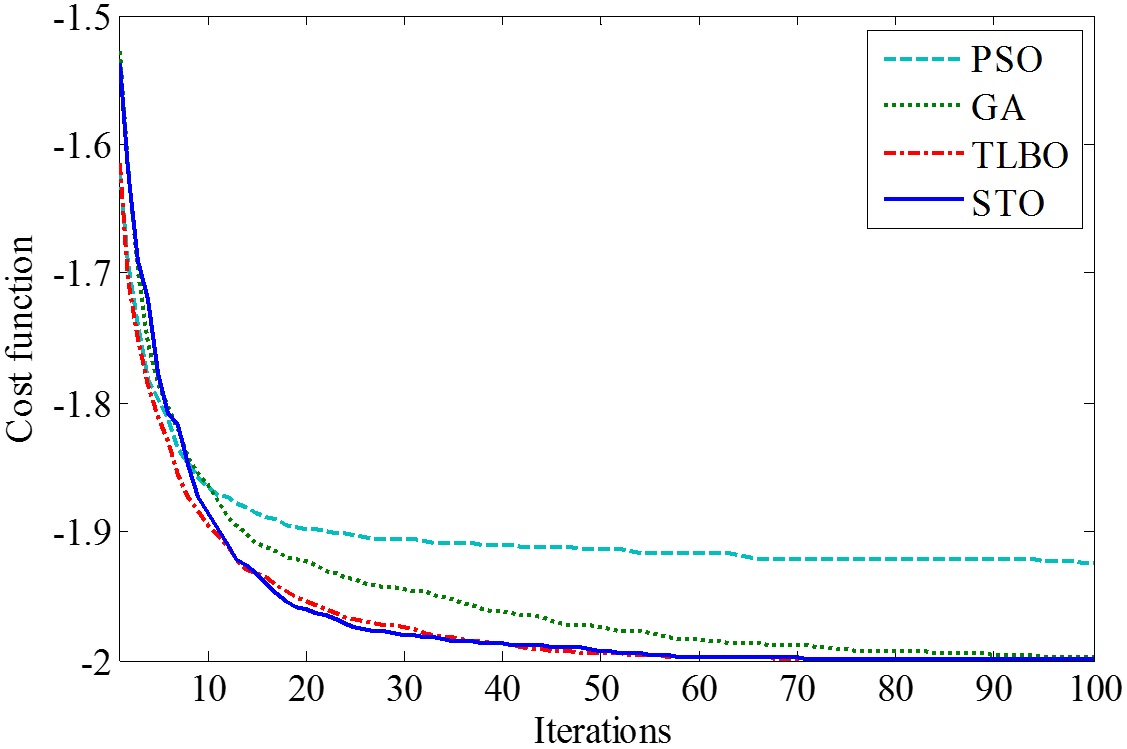}
  \caption{Convergence curves for \emph{Ripple25} function}
  \label{rip}
\end{figure}

\begin{figure}[h!]
 \centering
  \includegraphics[scale=0.3125]{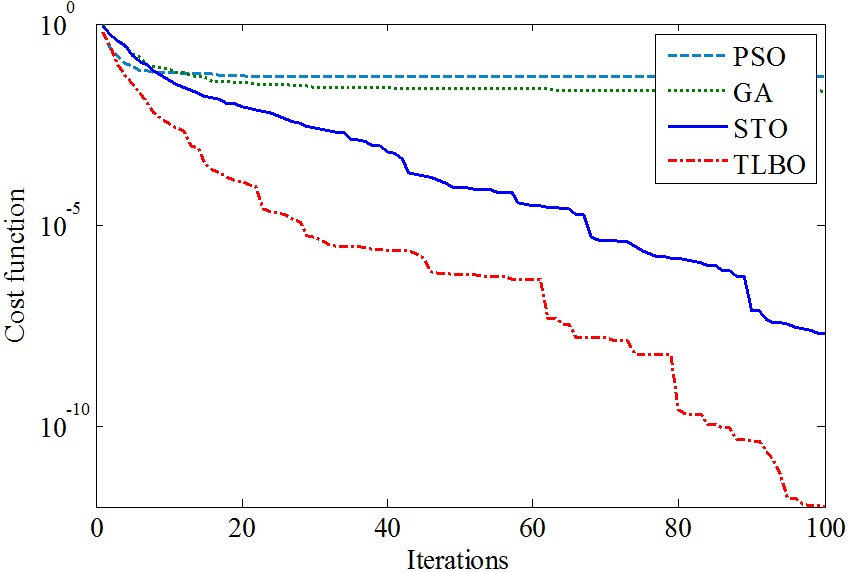}
  \caption{Convergence curves for \emph{Beale} Function}
  \label{kua2}
\end{figure}

Fig. \ref{kual1} to Fig. \ref{sty} show a qualitative comparison of the algorithms on optimization of \emph{EggHolder}, \emph{Ripple}, \emph{Beale}, and \emph{Modified Rosenbrock} functions, respectively. The convergence curves were obtained by averaging over 50 independent runs. As shown in the first three figures, the proposed algorithm has comparable performance with the TLBO as it was expected from the quantitative comparisons. In addition, STO in both tuned and randomized diameter versions, has best performance in the optimization of \emph{Modified Rosenbrock} (Fig. \ref{sty}).

\begin{figure}[h!]
 \centering
  \includegraphics[scale=0.23]{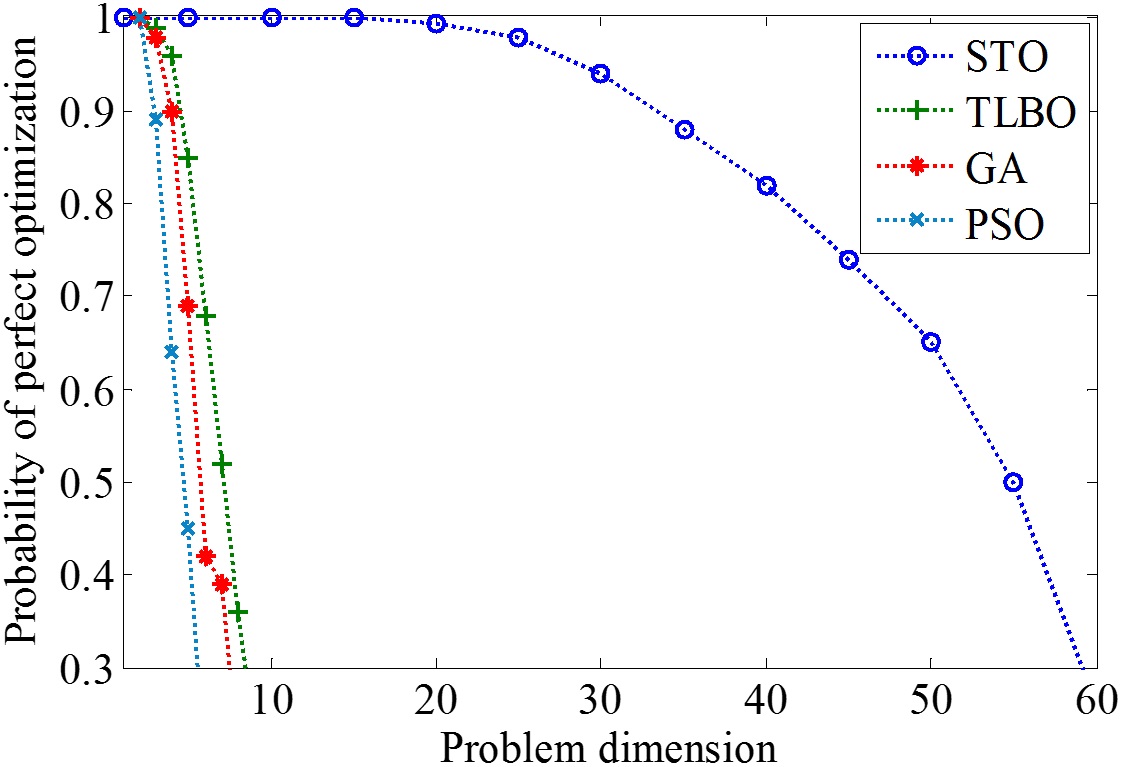}
  \caption{probability of perfect optimization of \emph{Styblinski-Tang} function in high dimensions}
  \label{highD}
\end{figure}

Fig. \ref{highD} shows the percentage of perfect optimizations (distortion less than $10^{-4}$) for the \emph{Styblinski-Tang} function respect to the dimension of problem. In two dimensions, this function consists of 3 local minimums and 1 global optimum (see Fig. \ref{Rast2}). As shown in Fig. \ref{highD}, all of the algorithms are completely successful in finding the global minimum at 2 dimensions. However, as the dimension increases, the PSO, GA, and TLBO were trapped in the local minimums such that their probability of perfect optimization dramatically decreases, while, the STO survives to explore new solutions in significantly higher dimensions. Maximum number of the allowed iterations for all of the algorithms was 5000. It was a limiting factor for the proposed algorithm, since it rarely fall into the local minimums of this problem and a better performance was possible by more iterations.

\begin{figure}[h!]
 \centering
  \includegraphics[scale=0.30]{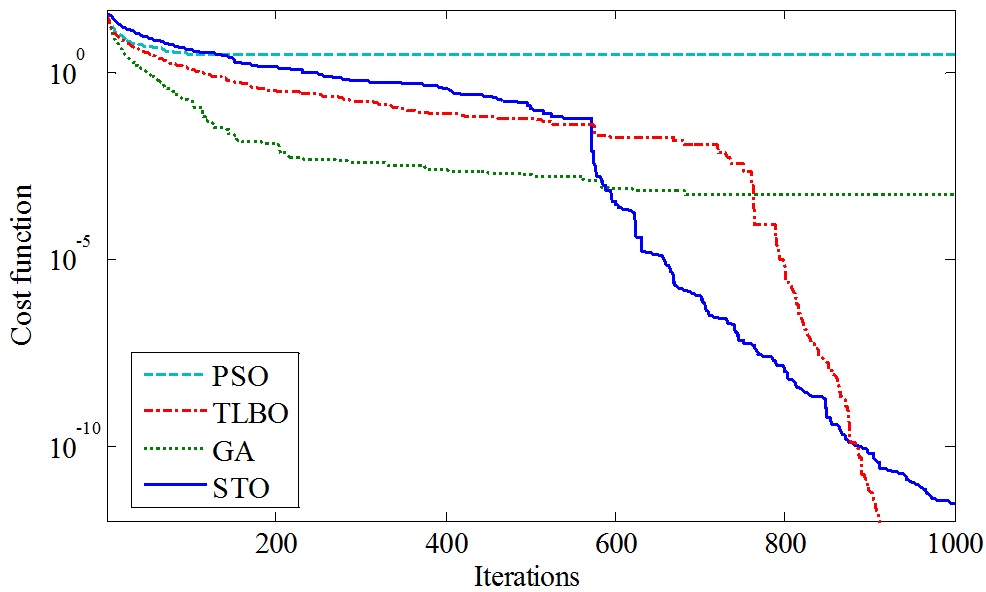}
  \caption{Convergence curves for \emph{Rastrigin} function in 5 dimensions}
  \label{rast}
\end{figure}

\begin{figure}[h!]
 \centering
  \includegraphics[scale=0.26]{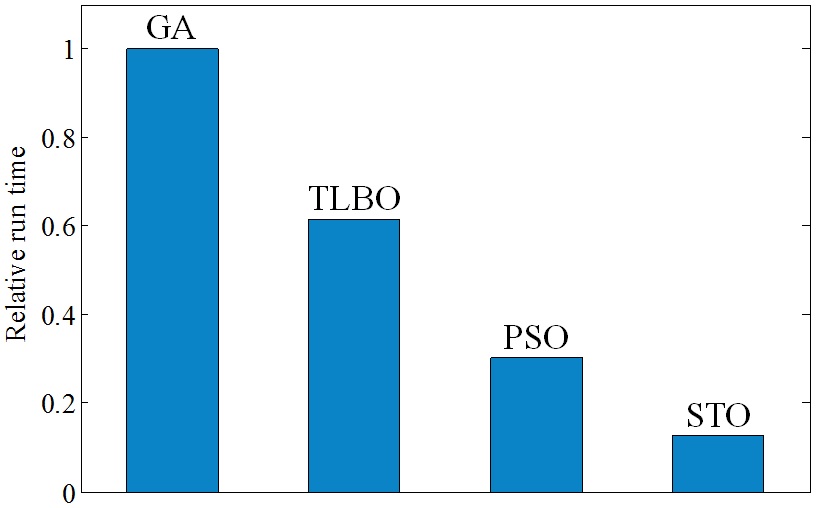}
  \caption{Comparison of relative run-time of the algorithms on last experiment}
  \label{run}
\end{figure}

In the other experiment, convergence curve of the STO in optimization of \emph{Rastrigin} function with 5 dimensions was evaluated in Fig. \ref{rast}. This function has a global minimum in the origin, with the cost function vlaue of zero. The results indicate a comparable performance of the proposed algorithm with the TLBO, both in quality and quantity (see Table \ref{tab}). In order to have a sense of the low computational complexity of the proposed algorithm, a normalized run-time of the algorithms for last experiment was provided at Fig. \ref{run}.

\section{Conclusion}

Synchronized air currents in a tornado leads to a significant reduction in the ambient temperature. It is a desired paradigm in optimization. Simulated tornado optimization (STO) mimicked two main air currents of a tornado in movement toward colder positions (better solutions). Participation of each air particle in only one current, or equivalently movement in just one search direction, led to low computational complexity in the proposed algorithm. Hence, it would be useful in the applications with time constraint. A simple self-adjustment is possible by randomization of the number of particles appointed to each current (randomization of the tornado diameter). Numerical results indicated the efficiency of the proposed algorithm in facing with different problems. As a future direction of research, more improvements can be achieved by adaptive reduction in the range of variations of the diameter. Also, a more realistic step size mechanism or turbulence model would be of great interest.\\

\appendix

\underline{Problems}:
\center
 (a) \emph{EggHolder} Function:

  \begin{displaymath}
    f=-(x_1+47)\sin\sqrt{|x_2+\frac{x_1}{2}+47|}-x_1\sin\sqrt{x_1-(x_2+47)}
  \end{displaymath}
  \begin{displaymath}
 -512\leq x_i \leq 512,~ f_{min}(512,404.2319)\approx-959.64
  \end{displaymath}

\center------------------------------------------------

  (b) \emph{Ripple25} Function:

  \begin{displaymath}
    f=\sum_{i=1}^{2}-e^{(-2\ln2(\frac{x_i-0.1}{0.8})^2)}\sin^6(5\pi x_i)
  \end{displaymath}
  \begin{displaymath}
 0\leq x_i \leq 1,~ f_{min}(0.1,0.1)=-2
  \end{displaymath}

\begin{center}
  ------------------------------------------------
\end{center}

  (c) \emph{Beale} Function:

  \begin{displaymath}
    f=(1.5-x_1+x_1x_2)^2+(2.25-x_1+x_1x_2^2)^2
      \end{displaymath}
      \begin{displaymath}
     +(2.625-x_1+x_1x_2^3)^2
      \end{displaymath}
  \begin{displaymath}
 -4.5\leq x_i \leq 4.5,~  f_{min}(3,0.5)=0
  \end{displaymath}

\center------------------------------------------------

(d) \emph{Modified Rosenbrock} Function:

  \begin{displaymath}
     f=74+100(x_2-x_1^2)^2+(1-x_1)-400e^{\frac{-(x_1+1)^2+(x_2+1)^2}{0.1}}
      \end{displaymath}

  \begin{displaymath}
-2\leq x_i \leq 2,~ f_{min}(-0.9,-0.95)=34.37
  \end{displaymath}

\center------------------------------------------------

(e) \emph{Styblinski-Tang} Function with $n$ decision variables:

  \begin{displaymath}
     f=\frac{1}{2}\sum_{i=1}^{n}(x_i^4-16x_i^2+5x_i),~ -5\leq x_i \leq 5
      \end{displaymath}
  \begin{displaymath}
f_{min}(-2.903534,\ldots,-2.903534)=-39.1661657037n
  \end{displaymath}

\center------------------------------------------------

(f) \emph{Rastrigin} Function with $n$ decision variables:

  \begin{displaymath}
     f=\frac{1}{2}\sum_{i=1}^{n}(x_i^2-10cos(2\pi x_i)+10)
      \end{displaymath}
  \begin{displaymath}
 -5.12\leq x_i \leq 5.12, ~ f_{min}(0,\ldots,0)=0
  \end{displaymath}
    %  \begin{displaymath}

 % \end{displaymath}\\

\center------------------------------------------------

\bibliographystyle{ieeetr}
\bibliography{STO}

\end{document}